\newcommand{\calN}{\mathcal{N}}
\newcommand{\calO}{\mathcal{O}}
\newcommand{\calP}{\mathcal{P}}

\newcommand{\calS}{\mathcal{S}}

\newcommand{\calW}{\mathcal{W}}

\newcommand{\bbP}{\mathbb{P}}

\def\SL{{\textrm{SL}}}

\documentclass[11pt]{amsart}
\usepackage{amscd, amsmath, amsthm, amssymb}

\theoremstyle{definition}     % italic or bold etc.

\theoremstyle{remark}

\numberwithin{equation}{section}

\begin{document}

\title[Luigi Cremona and cubic surfaces]
 {Luigi Cremona and cubic surfaces}

\author[Igor V.  Dolgachev]{Igor V. Dolgachev}
\address{Department of Mathematics, University of Michigan, Ann Arbor, MI 48109, USA}
\email{idolga@umich.edu}
\thanks{Research  is partially supported by NSF grant
DMS-0245203}

%\subjclass{[2000] 14J26}

\begin{abstract}
We discuss the contribution of  Luigi Cremona to the early development of the theory of cubic surfaces.
\end{abstract}

\maketitle

%%%%%%%%%%%%%%%%%%%%%%%%%%%%%%%%%%%%%%%%%%%

%\setcounter{section}{0}
\section{A brief history}

In 1911 Archibald Henderson wrote in his book \cite{Hen}

 ``While it is doubteless true that the classification of cubic surfaces is complete, the number of papers dealing with these surfaces which continue to appear from year to year furnish abundant proof of the fact that they still possess much the same fascination as they did in the days of their discovery of the twenty-seven lines upon the cubic surface.''

It is amazing that a similar statement can be repeated  almost a hundred years later.   Searching in MathSciNet for ``cubic surfaces'' and their close cousins ``Del Pezzo surfaces'' reveals 69 and 80 papers published in recent  10 years.  

Here are some of the highlights in the history of classical research on cubic surfaces before the work of Cremona. A good  source is Pascal's Repertorium  \cite{Pascal}.   

{\bf 1849}: Arthur Cayley communicates to George Salmon  that  a general cubic surface contains a finite number of lines. Salmon proves that the number of  lines must be equal to 27. Salmon's proof is  presented  in Cayley's paper \cite{Cayley}. In the same paper Cayley shows that a general cubic surface admits  45 tritangent planes, i.e.  planes planes which intersect the surface along the union of three lines. He gives a certain 4-parameter family of cubic surfaces for which the equations of tritangent planes can be explicitly found and their coefficients are rational functions in parameters.  

In a paper published in the same year and in the same journal \cite{Salmon1}   Salmon proves that not only a general but any nonsingular surface contains exactly 27 lines.  He also finds the number of lines on 11 different types of  singular surfaces. 

The discovery of 27 lines on a general cubic surface can be considered as the first non-trivial result on algebraic surfaces of order higher than 2. In fact, it can be considered as the beginning of modern algebraic geometry.

{\bf 1851}: John Sylvester claims without proof  that a general cubic surface can be written uniquely as a sum of 5 cubes of linear forms \cite{Sylvester}:
$$F_3 = L_1^3+L_2^3+L_3^3+L_4^3+L_5^3.$$
 This was proven ten years later  by Alfred  Clebsch \cite{Clebsch3}. The union of planes $L_i = 0$ will be known as the \emph{Sylvester pentahedron} of the cubic surface.

{\bf 1854}: Ludwig  Schl\"afli  finds about the 27 lines from correspondence with Jacob Steiner \cite{Graf}. In his letters he communicates  to Steiner some of his results on cubic surfaces which were published later in 1858 \cite{Schlafli1}. For example he shows that a general cubic surface has 36 double-sixes of lines. A \emph{double-six} is a pair of sets of 6 skew lines (\emph{sixes}) such that each line from one set is skew to a unique line from the other set making a bijection between the two sets. 
He introduces a new notation for the 27 lines ($a_i,b_i,c_{ij}, 1\le i <j\le 6$).  In these notations the double-sixes are 
\begin{equation}\notag
\small{\begin{pmatrix}a_1&\ldots&a_6\\
b_1&\ldots&b_6\end{pmatrix}},  \small{\begin{pmatrix}a_i&b_i&c_{kl}&c_{km}&c_{kn}&c_{kp}\\
a_k&b_k&c_{il}&c_{im}&c_{in}&c_{ip}\end{pmatrix}}, \small{\begin{pmatrix}a_i&a_k&a_l&c_{mn}&c_{mp}&c_{np}\\
c_{kl}&c_{ik}&c_{il}&b_p&b_{n}&b_m
\end{pmatrix}}\end{equation}

\smallskip
{\bf  1855}: Hermann Grassmann
proves  that three collinear nets of planes generate a cubic surface. More precisely, let 
$\calN_1,\calN_2,\calN_2$ be three general two-dimensional linear systems of planes in $\bbP^3$. Choose an isomorphism $\phi_i:\bbP^2\to \calN_i$ for each $i = 1,2,3$ and consider rational  map $\Phi:\bbP^2-\to \bbP^3$ defined by 
$$\Phi(x) = \phi_1(x)\cap \phi_1(x)\cap \phi_1(x).$$
Grassmann proves that its image is a cubic surface in $\bbP^3$. 

\smallskip
{\bf  1856}: Jacob Steiner \cite{Steiner}  introduces 120  subsets of 9 lines which form two pairs of triples of tritangent planes  (a \emph{ Triederpaar}). Each triple (a \emph{Trieder}) consists of tritangent planes intersecting along a line not on the surface. In Schl\"afli's notation they are 
$$\   \begin{matrix}a_i&b_j&c_{ij}\\
b_k &c_{jk}&a_j\\
c_{ik}&a_k&b_i
\end{matrix}, \qquad  \  \begin{matrix}c_{ij}&c_{kl}&c_{mn}\\
c_{ln} &c_{im}&c_{jk}\\
c_{km}&c_{jn}&c_{il}
\end{matrix},\qquad \ \begin{matrix}a_i&b_j&c_{ij}\\
b_k &a_l&c_{kl}\\
c_{ik}&c_{jl}&c_{mn}
\end{matrix}
$$
  This gives 120 essentially different representation of a general cubic surface by an equation of the form
$$L_1L_2L_3+L_4L_5L_6 = 0,$$
where $L_i$'s are linear forms. The existence of such a representation was already pointed out by Cayley \cite{Cayley} and Salmon \cite{Salmon1}. In the same memoir Steiner relates the Sylvester pentahedron with the linear system of polar quadrics  of the cubic surface. He also studies pencils of conics on a general cubic surface. Most of his results were stated without proof. 

\smallskip
{\bf 1860} Salmon \cite{Salmon2} and Clebsch \cite{Clebsch1},\cite{Clebsch2} find invariants of cubic forms in 4 variables. The ring of invariants is generated by invariants of degrees $8,16,24,32,40,100$ where the square of the invariant of degree 100 is a polynomial in the remaining invariants. In modern terms this implies that that the GIT-quotient of the projective space of cubic surfaces modulo the group of projective transformations is isomorphic to the weighted projective space $\bbP(1,2,3,4,5).$ 

\smallskip
{\bf 1862}: Fridericus  August \cite{August} proves that a general cubic surface can be projectively generated by three pencils of planes.
\smallskip 

\smallskip
{\bf 1863}: Ludwig Schl\"afli \cite{Schlafli2} classifies possible types of isolated singular points on cubic surfaces and their reality.

\smallskip
{\bf 1863}: Heinrich Schr\"oter \cite{Schroter} shows how to obtain the 27 lines on a cubic surface  obtained by Grassmann's construction and 
shows that there are exactly 6 triples of planes which intersect along a common line. The six lines form a half of  a double-six. In modern terms this implies that the cubic surfaces are isomorphic to the blow-up of six points in projective plane.

\smallskip
{\bf 1866} Alfred Clebsch \cite{Clebsch4} proves that a general cubic surface is the image of a birational map from projective plane given by the linear system of cubics through 6 points. Using this he shows that in the lines $a_i$ in Schl\"afli's notation are the images of the exceptional curves, the lines $b_i$ are the images of conics through 5 points, and the lines $c_{ij}$ are the images of lines through two points.
  
{\bf 1866}: Luigi Cremona and Rudolf Sturm are awarded a Prize of Berlin Academie established by Steiner for their work on cubic surfaces. A large part of their work was dedicated to supplying proof of the results announced by Steiner.

\smallskip
{\bf 1867} The work of  Sturm has been published \cite{Sturm1}.

\smallskip
{\bf 1868} The work of Cremona has been published \cite{Cremona1}.

 \section{Projective generatedness}
Let $X$ be a cubic surface with 3 skew lines $l_1,l_2,l_3$. Let $\calP_i, i= 1,2,3$ be three pencils of planes through the lines. A general plane $\pi_1\in \calP_1$ and a general plane $\pi_2\in \calP_2$ intersect  along a line $l$. The  line $l$ intersects the cubic surface at 3 points, one on $l_1$, one on $l_2$ and the third one  is the intersection point of two residual conics. Let $p$ be the third point and $\pi_3\in \calP_3$ be the plane through $p$. Then $\pi_1\cap \pi_1\cap \pi_3 = \{p\}$. Consider the rational map
$$\Psi:\bbP^1\times \bbP^1\times \bbP^1 -\to \bbP^3, (\pi_1,\pi_2,\pi_3)\mapsto \pi_1\cap \pi_1\cap \pi_3.$$
Let $R$ be the graph of the map $(\pi_1,\pi_2)\mapsto \pi_3$, where $\pi_3$ is the plane passing through $p = \pi_1\cap \pi_2\cap X$. Then the image of  $R$ under $\Psi$ is $X$.
The subvariety $R\subset \bbP^1\times \bbP^1\times \bbP^1$  defines a relation between 3 planes  which was called by  August \cite{August}  the {\it duplo-projective relation}.

Cremona wants to show that August's generatedness of $X$ implies the Grassmann generatedness by three nets of planes. This proves that every general cubic can be obtained by Grassmann's construction.   Here is his beautiful argument.   First, using the tritangent planes through the lines $l_1,l_2,l_3$, he finds the lines $l_{ij}$ on $X$ which intersects $l_i$ and $l_j$ but not three of them.   Then he takes a general plane $\pi$ and a triangle of lines $a_1,a_2,a_3$ on it. He chooses a projective isomorphism $a_i$ with $\calP_i$ such that the vertex $p_{ij} = a_i\cap a_j$ of the triangle   corresponds to the plane in $\calP_i$  containing the line $l_{ij}$. This does not yet fix the projective isomorphism. So, he chooses a point $\lambda_i\in a_i$ such that the planes $\pi_i$ from $\calP_i$ intersecting $a_i$ at $\lambda_i$ are in the correspondence $R$. Let $x_0 = \pi_1\cap \pi_2\cap \pi_3$. Now the projective isomorphisms $a_i\leftrightarrow \calP_i$ being fixed, he takes any point $x\in \pi$ and consider  the line joining $x$ with  the vertex $p_{ij}$ of the triangle. Each of three lines obtained in this way  intersects the opposite side $a_k$ at a point $\lambda_k(x)$ and hence defines a plane $\pi_k(x)\in \calP_k$. Let $x' $ be the intersection point of these planes.  He asks (\cite{Cremona1}, p. 70) {\it ``Quel est le lieu du point x?''} 
and proceeds to prove that it is equal to $X$. For this he shows that the locus $X'$ is a cubic surface containing 9 lines obtained from pairwise intersections of 6 planes spanned by  the lines $l_i$ and $l_{ij}$. Then he tries to show that the surface $X$ and $X'$ has also the point $x_0$ in common not lying on the union of the nine lines. This of course proves that $X = X'$. He uses that the  line $l_i$ correspond to the point $x$ taken on the side $a_i$ of the triangle and  the lines $l_{ij}$ correspond to the vertices of the triangle. But the argument for the latter is obscure. It was noticed by Sturm in his obituary of Cremona \cite{Sturm2}. He writes that he consulted on this point with Teodor Reye, who also did not understand the argument. The gap was fixed much later by Corrado Segre \cite{Segre2}. 
His proof is algebraic and very simple. One can write three plane-pencils in a duplo-projectuive relation as follows
$$\lambda_1 A_1+\lambda_2 A_2 = 0,  \mu_1 B_1+\mu_2 B_2 = 0, \gamma_1 C_1+\gamma_2 C_2 = 0$$
 $$\lambda_1\lambda_2\lambda_3 = \mu_1\mu_2\mu_3.$$
The equation of the cubic surface becomes 
$$F = A_1B_1C_1+A_2B_2C_2 = \det\begin{pmatrix}
     0 &A_1&A_2    \\
  B_2&0 & B_1\\
C_1&C_2&0 
\end{pmatrix} = 0$$
In general, if we have a $3\times 3$-matrix $M = (L_{ij})$ with  linear forms in 4 variables $t_0,t_1,t_2,t_3$ as its entries, the system of linear equations
$$L_1(\lambda) = \lambda_1L_{11}+\lambda_2L_{12}+\lambda_3L_{13} = 0$$
$$L_2(\lambda) = \lambda_1L_{21}+\lambda_2L_{22}+\lambda_3L_{23} = 0$$
$$L_2(\lambda) = \lambda_1L_{31}+\lambda_2L_{32}+\lambda_3L_{33} = 0$$
can be viewed as 3 projectively equivalent nets $V(L_\lambda)$ of planes in $\bbP^3$, a common solution is expressed by the determinant of the matrix $L = (L_{ij})$.
Rewrite this system of equations by collecting terms at the unknowns $t_i$:
$$t_0A_{11}+t_1A_{12}+t_2A_{13} +t_3A_{14}= 0$$
$$t_0A_{21}+t_1A_{22}+t_2A_{23} +t_3A_{24} =  0$$
$$t_0A_{31}+t_1A_{32}+t_1A_{33} +t_3A_{34} = 0.$$
Now the coefficients $A_{ij}$ are linear forms in $\lambda_1,\lambda_2,\lambda_3$. For any  $(\lambda_1,\lambda_2,\lambda_3)$ satisfying $F = \det L = 0$, a solution of the rewritten system of equations defines a point in $\bbP^3$. This shows that any point on the cubic surface $F = 0$ is obtained as the intersection of three planes each from its net of planes. This gives   Grassmann's representation of a cubic surface. Segre explains that Cremona's geometric argument is equivalent to this algebraic trick of rewriting the system of linear equations.  Segre also points out that the determinant equation $F = \det A$ exists even for singular cubic surfaces except when it acquires a special singular point, a double rational point of type $E_6$ in modern terminology.

Cremona derives other nice consequences of his geometric construction of projective degeneration of a cubic surface. First he considers a plane in  $\bbP^3 =\bbP(V)$ as a point in the dual space $\check{\bbP}^3 = \bbP(V^*)$. Intersection of three planes defines a rational map $f$ from   $(\check{\bbP}^3)^3$ to $\bbP^3$. Now choose a projective generation of a cubic surface $X$ by three projectively equivalent nets of planes $V(L_i\lambda)$ as above. Assigning to a point $x\in \bbP^3$ a triple of planes $\Pi_i$ from each net defines a  rational map $\phi$  from $\bbP^3$ to $(\check{\bbP}^3)^3$. Composing it with the map $f$, we get a rational map
$$T:\bbP^3 -\to \bbP^3.$$
It is the famous  {\it cubo-cubic Cremona transformation} of $\bbP^3$. It is easy to see that $\phi$ is defined by a linear map $V\to V^*\otimes V^*\otimes V^*$ and the map $f$ is defined by a linear map $V^*\otimes V^*\otimes V^*$ which factors through the symmetric product $S^3(V^*)$.  This shows that the map $T$ is given by a 3-dimensional linear system  of cubics. Its base locus is the pre-image of the locus of linear dependent triples of planes. It is a curve $C$ of degree 6 of genus 3. The image of a plane is a surface of degree 3. The map is not defined on the intersection of $C$ with the plane. This is a set of 6 points.
This recovers Clebsch's result that  the rational map $\Phi$ is given by a linear system of cubics through 6 points in the plane. 

Using this representation Cremona easily lists all 27 lines 
{\it ``Voil\`a donc les 27 droites de la surface $F_3$''} \cite{Cremona1}, p. 67.

This allowed him to reconstruct Steiner's Triederpaare, tritangent planes, and double-sixes. He makes an important observation that a choice of a double-six determines two 3-dimensional linear systems of curves of degree 6 is residual to each other with respect to complete intersections of $X$ with a quadric.  As we know now such a linear system defines a determinantal representation of the cubic equation.

Steiner announces, and Cremona and Sturm prove, that each Triederpaar of tritangent planes defines two more Triederpaare such that together they consist of the 27 lines. The number of such triples is equal to 40.  Each triple defines 9 tritangent planes which together contain the 27 lines. The sets of 9 tritangent planes with this property were first considered by Camille Jordan \cite{Jordan1}. He called such a set a \emph{enneaedro}.  In \cite{Cremona3} Cremona divides all enneaedra in two kinds. An enneaedro of the first kind is obtained as above from a triple of Triederpaare.  There are 40 of them.  Only those were found by Jordan. An enneaedro of the second kind i can be uniquely divided into the union of three triples of planes, each triple belongs to a unique triple of Triederpaare as above.

\section{Cremona's hexaedral equations}
The last paper of Cremona devoted to cubic surfaces is his Math. Ann. paper of 1878 \cite{Cremona4}. It was shown by Reye \cite{Reye1}  that the equation of a general cubic surface $F_3$ can be expressed in $\infty^4$ ways as a sum of 6 linear forms:
$$F_3 = \sum_{i=1}^6 L_i^3 = 0$$
(a {\it hexaedral equation}). The linear forms $l_i$ define a map $\bbP^3\to \bbP^5$ with the image a linear susbspace defined by 2 linear independent relations between the linear forms:
$$\sum_{i=0}^5 a_ix_i = 0, \quad \sum_{i=0}^5 b_ix_i = 0.$$
The surface $F_3$ is mapped to the linear section of the cubic hypersurface
$$\sum_{i=0}^5 x_i^3 = 0.$$
Cremona  shows that among $\infty^4$ representations of $F$ one can always choose a representation such that one of the equations of the linear space is 
$$x_0+\ldots+x_5 = 0.$$
In fact it can be done in 36 essentially different ways corresponding to a choice of a double-six on $F_3$. This is done as follows. A choice of a Steiner's Triededrpaar of tritangent planes gives a Cayley-Salmon equation of $X$ of the form  
$F = PQR+STU$. Then one looks for constants such that after scaling the linear forms they add up to zero. Write
$$P' = pP, Q' = qQ, R' = rR, S' = sS, T' = tT, U' = uU.$$
Since $V(F)$ is not a cone, four of the linear forms are linearly independent. After reordering the linear forms, we may assume that the forms $P,Q,R,S$ are linearly independent. Let 
$$T = aP+bQ+cR+dS, \ U = a'P+b'Q+c'R+d'S.$$
The constants $p,q,r,s,t,u$ must satisfy  the following system of equations
\begin{eqnarray*}
p+ta+ua' &=& 0\\
q+tb+ub' &=& 0\\
r+tc+uc'&=&0\\
s+td+ud'&=&0\\
pqr+stu &=&0.
\end{eqnarray*}
The first 4 linear equations allow us to express linearly all unknowns in terms of two, say $t,u$. Plugging in the last equation, we get  a cubic equation in $t/u$. Solving it, we get a solution. Now set
$$x_1 = Q'+R'-P', \quad x_2 = R'+P'-Q',\quad x_3 = P'+Q'-R',$$
$$x_4 = T'+U'-S',\quad x_5 = U'+S'-T', \quad x_6 = S'+T'-U'.$$
One checks that these six linear forms satisfy the equations 
$$x_1^3+\ldots+x_6^3 = 0, \quad x_1+\ldots+x_6 = 0.$$
   A {\it Cremona hexaedral equation} defines a set of 15 lines given by the equations 
$$x_i+x_j = x_k+x_l = x_m+x_n = 0, \ \sum_{i=1}a_ix_i = 0.$$
One check that the complementary set is a double-six. 
If two surfaces given by hexahedral  equations define the same double-six, then they have common 15 lines. Obviously this is impossible. Thus the number of different hexahedral equations of $X$ is less or equal than 36. Now consider the identity
$$(x_1+\ldots+x_6)\bigl((x_1+x_2+x_3)^2+(x_4+x_5+x_6)^2-(x_1+x_2+x_3)(x_4+x_5+x_6)\big) $$
$$=(x_1+x_2+x_3)^3+(x_4+x_5+x_6)^3 = x_1^3+\ldots+x_6^3$$
$$+3(x_2+x_3)(x_1+x_3)(x_1+x_2)+3(x_4+x_5)(x_5+x_6)(x_4+x_6).$$
It shows that  Cremona hexahedral  equations define a Cayley-Salmon equation
$$(x_2+x_3)(x_1+x_3)(x_1+x_2)+(x_4+x_5)(x_5+x_6)(x_4+x_6) = 0,$$
where we have to eliminate one unknown with help of the equation $\sum a_ix_i = 0$.
Applying permutations of $x_1,\ldots,x_6$, we get 10 Cayley-Salmon equations of $S$. Each 9 lines formed by the corresponding Steiner's Triederpaar  are among the 15 lines determined by the hexahedral  equation. It follows from the classification of Steiner's Triederpaare that we have 10 such pairs composed of lines $c_{ij}$'s. Thus a choice of  Cremona hexahedral  equations defines exactly 10 Cayley-Salmon equations of $S$. Conversely, it follows from Cremona's proof  from above that  each Cayley-Salmon equation gives 3 Cremona hexahedral  equations (unless the cubic equation has a multiple root). Since we have 120 Cayley-Salmon equations for $S$ we get $36 = 360/10$ hexahedral  equations for $S$. They match with 36 double-sixes.

The cubic 3-fold $\calS_3$ in $\bbP^5$ given by equations 
$$\sum_{i=0}^5x_i^3 = \sum_{i=0}^5x_i = 0$$
is the Segre cubic primal \cite{Segre1}.  It has 10 nodes, the maximum possible for a cubic 3-fold.
In 1915 Arthur Coble computes invariants of 6 points on $\bbP^1$ and proves (in modern terms) that  $\calS_3$ is isomorphic to the GIT-quotient $(\bbP^1)^6/\!/\SL(2)$ \cite{Coble}. A choice of a six skew lines defines a representation of $F_3$ as the blow-up of 6 points 
$p_1,\ldots,p_6$ in $\bbP^2$. We order them. Consider the projection of these 6 points to $\bbP^1$ from a variable point $x\in \bbP^2$. We get a map from $\bbP^2$ to $S_3$. Coble proves that the image is a hyperplane section $\sum_{i=0}^5 a_ix_i = 0.$ This hyperplane depends only on the projective equivalence class of the six points. Also if we replace the 6 points by the six points coming from 6 skew lines which form with the previous set of skew lines a double six, we get the same hyperplane section. The moduli space of nonsingular cubic surfaces together with a choice of an ordered set of skew lines  is isomorphic to the open subset of the GIT-quotient $(\bbP^2)^6/\!/\SL(3)$. The involution which interchanges the 6 lines with the dual 6 lines, has the quotient isomorphic to $\bbP^4$. This is the $\bbP^4$ formed by the coefficients $(a_0,\ldots,a_5)$ in Cremona hexaedral equations!

\section {Desmic quartics}
 Fix a tritangent plane $\pi$ on $X$ which consists of 3 lines $l_1,l_2,l_3$. A plane $\pi_i$ through $l_i$ has a residual conic $C_i$. Since the reducible curve
$l_1+l_2+l_3$ is cut out by a plane, $C_1+C_2+C_3$ is cut out by a quadric. Thus we have a map 
$$\calP_1\times \calP_2\times \calP_3 \to |\calO_{\bbP^3}(2)|,
$$
where $\calP_i$ is the pencil of planes through the line $l_i$. Since the map does not depend on the order of lines, it factors through a linear map $\bbP^3\to  |\calO_{\bbP^3}(2)|$. This defines a web $\calW$ of quadrics in $\bbP^3$.  Choosing the plane $\pi_i$ among the 4 tritangent planes through $l_i$ except the plane $\pi$,  gives a quadric intersecting $X$ along the union of 6 lines. By changing $\pi$, one  obtains 45 sets  of 48 nonsingular quadrics with this property. Each quadric belongs to 6 sets, altogether giving $360$ quadrics.  This beautiful result of Steiner was given without proof. The proofs were supplied by Cremona and Sturm.

They went further by proving that the Steinerian  surface of the web $\calW$   (the set of nodes of quadrics from the web) is a  quartic surface $K$ with 12 nodes. The twelve nodes lie on $X$ and are equal to the intersection points of 12 line-pairs corresponding to the twelve  tritangent planes from above. The twelve nodes of $K$ are the vertices of 3 desmic tetrahedra in $\bbP^3$ (\emph{desmic} means that their equations are linear dependent). In particular, the 135 intersection points of the 27 lines on $X$ can be grouped in 45 sets of 12 points which are the vertices of 3 desmic tetrahedra. It was proven by C. Jessop  \cite{Jessop} in 1916  that the \emph{desmic quartic} $K$ is birationally isomorphic to the Kummer surface associated to the product of two isomorphic elliptic curves.

\section{Pascal's  Hexagram}  
A well-known Pascal's Theorem from plane projective geometry asserts that the opposite sides of a hexagon inscribed in a conic $C$ intersect at three collinear points. In \cite{Cremona3} Cremona observes that this configuration can be obtained by projecting a cubic surface with an ordinary double point $O$ with center at $0$. The 6 lines $l_1$,\ldots,$l_6$ passing through $O$ are projected to six points $p_1,\ldots,p_6$ lying on a conic. Each of 15 pairs of lines $l_i,l_j$ defines a tritangent plane with the residual line $l_{ij}$. The remaining 15 tritangent planes formed by lines $l_{ij},l_{kl},l_{mn}$ where the index sets have no common elements. The lines 
$l_{ij}$ are projected to lines $\ell_{ij}$ through the points $p_i,p_j$. Two tritangent planes  not containing the point  $O$ and with no common lines define a hexagon inscribed in the conic $C$. The intersection line of the two tritangent planes is projected to the line joining the intersection points of opposite sides of the hexagon. In this way a nodal cubic surface $X$ defines 60 Pascal's Hexagrams. Next Cremona observes that all of this can be extended to a nonsingulat cubic surface. A choice of a double-six defines a Cremona hexagonal form of the equation of $X$. The remaining 15 lines can be indexed by the sets $\{ij,kl,mn\}$. They lie in the union of 15 tritangent planes $\Pi_{ij}$ which can be indexed by the sets $\{ij\}$. Two tritangent planes $\Pi_{ij}$ and $\Pi_{kl}$ form a \emph{Cremona pair} if they  intersect along a line not contained in $X$. This terminology belongs to Reye (\cite{Reye2}, p. 218) who also calls the intersection line a \emph{Pascal line}.  There are 60 Cremona pairs. Projecting from a general point on $X$ we obtain 60 hexagons with opposite sides intersecting at 3 points lying on a line (the projection of a Pascal line). So we get 60 Pascal's hexagrams although there is no conic in which the hexagon is inscribed! 

The 15 tritangent planes $\Pi_{ij}$ can be realized as the faces of 6 pentahedra $P_i$ (the \emph{Cremona pentahedra}), its 60 edges are the Pascal lines. All subsets of the set of 45 tritangent planes with the property that no two planes in the set intersect along a line on $X$ were found by Eugenio Bertini \cite{Bertini}.

As it turned out much later all this fascinating  combinatorics of the set of 27 lines has a nice group-theoretic interpretation.  It was proven by Jordan \cite{Jordan2} that the group of symmetries of the incidence graph  of the 27 lines is 
realized as the Galois group of the equation defining the 27 lines whose coefficients are rational functions in the coefficients  of the equation of a general cubic surface. In modern terms this group is isomorphic to the Weyl group $W(E_6)$ of  root system of type $E_6$. The double-sixes correspond to pairs of opposite roots, so that the stabilizer of such a pair is a maximal subgroup of index 36. The stabilizer of a line is a maximal subgroup of index 27. The stabilizer of a tritangent plane is a maximal subgroup of index 45. The stabilizer of a triple of Steiner's Triederpaare containing all 27 lines is a maximal subgroup of index 40. In this way all (except one of index 40) maximal subgroups get a geometric interpretation in terms of the geometry of 27 lines on a cubic surface.

\section{Real lines on a cubic surface}
The last chapter of Cremona's memoir is dedicated to the questions of reality of a cubic surface and its lines. Schl\"afli \cite{Schlafli2} and August \cite{August} have already distributed nonsingular cubic surfaces into  5 different species. For example, surfaces of the first kind have all 27 lines  real. Surfaces of the second kind have 15 lines are real and 6 complex conjugate pairs of lines form a double-six. Cremona and Sturm give geometric proofs of these results. Cremona's  classification takes into account also the reality of tritangent planes and double-sixes.
 
\begin{itemize}
\item[(1)] All is real (lines, tritangent planes, double-sixes);
\item[(2)] 15 lines and  15 tritangent planes are real. There are 15 double-sixes  which consist of real sixes, each six consists of 6 real lines and a pair of complex conjugate lines.  One double-six consists of two complex conjugate sixes.
\item[(3)] 7 lines and 5 tritangent planes are real. There are 6 double-sixes  which consist of two real sixes  formed by 2 real lines and two pairs of complex conjugate lines. There are also 2 double-sixes,  each consists of a pair of complex conjugate sixes. 
\item[(4)] 3 lines and 7 tritangent planes are real. There is a unique double-six which consists  of two real sixes  formed by three pairs of complex conjugate lines. There are also 3 real-double sixes,  each consists of a complex conjugate pair of sixes.
\item[(5)] 3 lines and 13 tritangent planes are real. No real double-sixes.
\end{itemize}

Cremona and Sturm prove that only surfaces of the fifth kind cannot be generated  by three nets of planes as in  Grassmann's construction. The reason is that on such a surface there are no sets of 6 skew lines defined over reals.  

The argument of Cremona is very geometric. He first observes that a cubic surface can be reconstructed from a Triederpaar of tritangent planes and a point by using the Cayley-Salmon equation. A real surface is obtained if we find a Triederpaar and a point defined over reals. A real tritangent plane is formed by either 3 real lines or one real line and a pair of complex conjugate lines. A real Triederpaar is defined by a pair of real tritangent planes intersecting 
along a line not contained in the surface. Cremona proves that it is always possible to find a real Triederpaar on a real surface. Next he proceeds to consider different cases.

Case 1: The nine  lines defined by a real Triederpaar are all real.

A Triederpaar contains 3 triples of skew lines. A triple of skew lines defines a unique quadric $Q$ containing them in its ruling. The quadric $Q$ intersect the surface along 3 lines in another ruling. If $Q$ can be found such that the new 3 lines are all real, we get  type (1) from above. If $Q$ defines 3 lines such that one of them is real and 2 are complex conjugate, we obtain case (2). 

Case 2: One  of the Trieders in the pair consists of real planes, the other one consists of a real plane and two complex conujgates. This case leads to types (4) and (5).

Case 3: Each Trieder consists of  one real plane and two complex conjugate planes. This leads to cases (3) and (4).


\begin{thebibliography}{AMR}
\bibitem[Aug]{August} F. August, \emph{Discusitiones de superfieciebus tertii ordinis} (in Latin), Diss. Berlin. 1862. Available on the web from the G\"ottingen Mathematical Collection.

\bibitem[Ber]{Bertini} E. Bertini, \emph{Contribuzione alla teoria delle 27 rette e dei 45 piani tritangenti di una superficie de $3^\circ$ ordini}, Ann. di Mat. (2), {\bf 12} (1884), 301--346.

\bibitem[Cay]{Cayley} A. Cayley, \emph{On the triple tangent planes of surfaces of the third order}, Cambridge and Dublin Math. J., {\bf 4} (1849), 118--138.


\bibitem[Cle1]{Clebsch1} A. Clebsch, \emph{Zur Theorie der algebraischer Fl\"achen}, Journ.  f\"ur  reine und angew. Math., {\bf 58} (1860), 93--108.

\bibitem[Cle2]{Clebsch2} A. Clebsch, \emph{Ueber eine Transformation der homogenen Funktionen dritter Ordnung mit vier Ver\"anderlichen}, Journ.  f\"ur  reine und angew. Math., {\bf 58} (1860), 109--126.

\bibitem[Cle3]{Clebsch3} A. Clebsch, \emph{Ueber die Knotenpunkte der Hesseschen Fl\"ache, insbesondere bei Oberfl\"achen dritter Ordnung}, Journ.  f\"ur reiner und angew. Math., {\bf 59} (1861), 193--228.

\bibitem[Cle4]{Clebsch4} A. Clebsch, \emph{Die Geometrie auf den Fl\"achen dritter Ordnung}, Journ.  f\"ur reine und angew. Math., {\bf 65}, (1866), 359--380.

\bibitem[Cob]{Coble} A. Coble, \emph{Point sets and allied Cremona groups}, Trans. AMS, {\bf 16} (1915), 155--198.

\bibitem[Cre1]{Cremona1} L. Cremona, \emph{M\'emoire de g\'eom\'etrie pure sur les surfaces du troisi\'eme ordre}, Journ.  des Math\'ematiques pures et appliqu\'ees,  {\bf 68} (1868), 1-- 133 (Opere matematiche di  Luigi Cremona, Milan, 1914, t. 3, pp.1-121) [German translation: \emph{Grunz\"uge einer allgeimeinen Theorie der Oberfl\"achen in synthetischer Behandlung}, Berlin, 1870].

\bibitem[Cre2]{Cremona2}  L. Cremona, \emph{Ventisette rette di una superficie del terze ordine}, Rendiconti Inst. Lombardo, (2), {\bf 3} (1870), 209--219 (Opere, t. 3, pp. 167--176).


\bibitem[Cre3]{Cremona3}  L. Cremona, \emph{Theoremi stereometrici dai quali si deducono le propriet\`a dell'esagrammo di Pascal}, Accad. Lincei, Memorie della classe di scienze mat. e nat. (3), {\bf 1} (1876-77), 854--874 (Opere, t. 3, pp. 406--429).

\bibitem[Cre4]{Cremona4}  L. Cremona, \emph{Ueber die Polar-Hexaeder bei den Fl\"achen dritter ordnung}, Math. Ann. {\bf 13} (1878), 301--304 (Opere, t. 3, pp. 430-433).

\bibitem[Graf]{Graf} J. H. Graf, \emph{Der Briefwechsel zwichen Jakob Steiner und Ludwig Schl\"afli}, Mitteilung der Naturforschung Gesselschaft. Bern. 1896.

\bibitem[Gras]{Grassmann} H. Grassmann, {\it Die stereometrische Gleichungen dritten Grades und die dadurch erzeugen Oberfl\"achen}, Journ.  f\"ur reiner und angew. Math., {\bf 49} (1856), 47-- 65.

\bibitem[Hen]{Hen} A. Henderson, \emph{The  twenty-seven lines upon the cubic surface}, Cambridge, 1911.

\bibitem[Jes]{Jessop} C. Jessop, \emph{Quartic surfaces with singular points}, Cambridge Univ. Press, 1916.

\bibitem[Jor]{Jordan1} C. Jordan, \emph{Sur une nouvelle combinaison des 27 droites d'une surface du troisime ordre}, Comp. Rendus Acad. Sci. Paris, {\bf 70} (1870), 326--328.

 \bibitem[Jor]{Jordan2} C. Jordan,\emph{Trait\'e des substitutions et des \'equations alg\'ebriques}. Paris. Gautier-Villars. 1870.


\bibitem[Pas]{Pascal} E. Pascal, \emph{Repertorium der H\"oheren Mathematik}, Bd. 2, Teubner, Leipzig und Berlin, 1922.

\bibitem[Rey1]{Reye1} T. Reye, \emph{Geometrische Beweis des Sylvesterschen Satzes:``Jede quatern\"are cubische Form is darstbellbar als Summe von f\"unf Cuben linearer Formen''}, Journ.  f\"ur reiner und angew. Math., {\bf 78} (1874), 114--122.

\bibitem[Rey2]{Reye2} T. Reye, \emph{Die Geometrie der Lage}, t. 3, Leipzig. 1886.

\bibitem[Sal1]{Salmon1} G. Salmon, \emph{On the triple tangent planes to a surface of the third order}, Cambridge and Dublin Math. J., {\bf 4} (1849), 252--260

\bibitem[Sal2]{Salmon2} G. Salmon, \emph{On quaternary cubics}, Phil. Trans. of Roy. Soc.  London, {\bf 150} (1860), 229--237.

\bibitem[Schl1]{Schlafli1} L. Schl\"afli, \emph{An attempt to determine the twenty-seven lines upon a surface of the third order and to divide such surfaces into species in reference to the reality of the lines upon the surface}, Quart. J. Math., {\bf 2} (1858), 55--65, 110--121. 

\bibitem[Schl2]{Schlafli2} L. Schl\"afli, \emph{On the distributionn of surfaces of the third order into species, in reference to the absence or presense of singular points, and the reality of their lines}, Phil. Trans. of Roy. Soc. London, {\bf 6} (1863), 201--241. 

\bibitem[Schr]{Schroter} H. Schr\"oter, \emph{Nachweis der 27 Geraden auf der allgemeinen Oberfl\"achen dritter Ordnung}, Journ.  f\"ur reiner und angew. Math., {\bf 62} (1863), 265--280.

\bibitem[Seg1]{Segre1} C. Segre, \emph{Sulla variet\`a cubica con dieci punti doppi dello spazio a quattro dimensioni}, Atti Aca. Scienz. Torino, {\bf 22} (1886/87), 791--801
(Opere, Edizioni Cremonese, Roma, 1963,  v. IV, pp. 88--98).

\bibitem[Seg2]{Segre2} C. Segre, \emph{Sur la g\'en\'eration projective des surfaces cubiques}, Archiv der Math. und Phys., (3) {\bf 10} (1906), 209--215 (Opere, v. IV, pp. 188--196).

\bibitem[Ste]{Steiner} J. Steiner, {\it Ueber die Fl\"achen dritten Grades}, Journ.  f\"ur reiner und angew. Math., {\bf 53} (1856), 133--141 (Gesammelte Werke, Chelsea, 1971, vol. II, pp. 649--659).






\bibitem[Stu1]{Sturm1} R. Sturm, {\it Synthetische Untersuchungen \"uber Fl\"achen dritter Ordnung}, Teubner, Leipzig, 1867.

\bibitem[Stu2]{Sturm2} R. Sturm, {\it Luigi Cremona}, Archiv der Math. und Phys., (3), {\bf 8} (1904), 195--213.

\bibitem[Syl]{Sylvester} J. Sylvester, \emph{Sketch of a memoir on elimination, transformation, and canonical forms}, Cambridge and Dublin Math. J., {\bf 6} (1851), 186--200.






\end{thebibliography}
\end{document}